# NOTES ON HILBERT-KUNZ MULTIPLICITY OF REES ALGEBRAS

Kazufumi Eto
Department of Mathematics
Nippon Institute of Technology
Miyashiro, Saitama 345–8501, Japan
e-mail: `etou@nit.ac.jp`

and

Ken-ichi Yoshida
Graduate School of Mathematics
Nagoya University
Chikusa-ku, Nagoya 464–8602, Japan
e-mail: `yoshida@math.nagoya-u.ac.jp`

ABSTRACT. In this paper, we estimate the Hilbert-Kunz multiplicity of the (extended) Rees algebras in terms of some invariants of the base ring. Also, we give an explicit formula for the Hilbert-Kunz multiplicities of Rees algebras over Veronese subrings.

## INTRODUCTION

Throughout this paper, let $(A, \mathfrak{m}, k)$ be a commutative Noetherian local ring with unique maximal ideal $\mathfrak{m}$ of characteristic $p > 0$ with $d := \dim A \geq 1$. This paper is devoted to studying the Hilbert-Kunz multiplicities of (extended) Rees algebras over $A$.

The notion of Hilbert-Kunz multiplicity (denoted by $e_{HK}(I)$) has been defined by Kunz [14,15] and was formulated by Monsky [17] explicitly.

In 1980's, Hochster and Huneke [10] have introduced the notion of tight closure and showed that the tight closure of an ideal is the largest ideal containing the ideal having the same Hilbert-Kunz multiplicity (under some mild conditions); see also [19] or [12, Theorem 5.3]. Furthermore, in [24], K.-i Watanabe and the second-named author have proved that an unmixed local ring whose Hilbert-Kunz multiplicity (denoted by $e_{HK}(A)$) is equal to one is regular; see also [18, (40.6)]. These facts indicate that there exist the parallels between the notion of integral closures and that of tight closures in terms of the comparison between the notion of multiplicity (denoted by $e(I)$) and that of Hilbert-Kunz multiplicity.

Typeset by $\mathcal{A}_{\mathcal{M}}\mathcal{S}$-TEX

Hilbert-Kunz multiplicity is a sort of "multiplicity", but it is *not* integer in general. Thus it is important to determine the range of the value of the Hilbert-Kunz multiplicity. For example, let $A$ be a hypersurface of multiplicity 2. Then $A$ is F-rational (resp. *not* F-rational) if and only if $1 < e_{HK}(A) < 2$ (resp. $e_{HK}(A) = 2$).

In this context, we consider the Hilbert-Kunz multiplicities of (extended) Rees algebras. We recall the notion of Rees algebras. The algebra $R(I) := \oplus_{n \in \mathbb{N}} I^n = A[It]$ (resp. $R'(I) := \oplus_{n \in \mathbb{Z}} I^n = A[It, t^{-1}]$) is called the *Rees algebra* (resp. the *extended Rees algebra*) of $A$ with respect to $I$ (or $\{I^n\}$). Several properties (e.g. Cohen–Macaulay, Gorenstein etc.) of these algebras have been investigated by many authors. In particular, as for multiplicity, the following fact is known.

**Fact.** *Let $I \subseteq A$ be an $\mathfrak{m}$-primary ideal. Put $G(I) := \oplus_{n \geq 0} I^n/I^{n+1}$, the associated graded ring of $I$ and $\mathfrak{M} := \mathfrak{m}R(I) + R(I)_+$. Also, we put $e(R(I)) := e(R(I)_{\mathfrak{M}})$ (resp. $e(G(I)) := e(G(I)_{\mathfrak{M}G(I)})$). Then*
   (1) $e(I) = e((t^{-1}, R'(I)_+)R'(I)) = e(G(I))$.
   (2) $e(R(\mathfrak{m})) = d \cdot e(A)$; *see e.g.* [22].

So it is natural to ask whether or not the similar formula holds for Hilbert-Kunz multiplicity. In fact, we propose the following question.

**Question.** *Let $e_{HK}(I)$ denote the Hilbert-Kunz multiplicity of $I$. Then*
   (1) *Do inequalities $e_{HK}(I) \leq e_{HK}((t^{-1}, R'(I)_+)R'(I)) \leq e_{HK}(G(I))$ always hold?*
   (2) *Do inequalities $e(A) \leq e_{HK}(R(\mathfrak{m})) \leq c(d) \cdot e(A)$ always hold, where $c(d)$ is a positive real number depending on $d = \dim A$ only?*
   (3) *When do equalities hold in (1) or (2)?*

In this paper, we prove the following two theorems as partial answers to the above question.

**Theorem 1.** *Let $(A, \mathfrak{m})$ be a local ring of characteristic $p > 0$ with $\dim A \geq 1$. Then for any $\mathfrak{m}$-primary ideal $I \subseteq A$, we have*

$$e_{HK}(I) \leq e_{HK}((t^{-1}, R'(I)_+)R'(I)) \leq e_{HK}(G(I)) \leq e(I).$$

**Theorem 2.** *Let $(A, \mathfrak{m})$ be a local ring of characteristic $p > 0$ with $d = \dim A \geq 1$. Then we have*

$$e_{HK}(R'(\mathfrak{m})) \leq e_{HK}(R(\mathfrak{m})) \leq c(d) \cdot e(A), \quad \text{where} \quad c(d) = d\left(\frac{1}{2} + \frac{1}{(d+1)!}\right).$$

*Further, $e_{HK}(R(\mathfrak{m})) = c(d) \cdot e(A)$ if and only if $e_{HK}(A) = e(A)$.*

By Theorem 1 and Theorem 2, we have that $e_{HK}(A) \leq e_{HK}(R'(\mathfrak{m})) \leq e_{HK}(R(\mathfrak{m}))$. Also, we will prove that $e(A) \leq e_{HK}(R(\mathfrak{m}))$ holds in case of two-dimensional Cohen–Macaulay local rings. However, this inequality does not hold in general if $\dim A \geq 3$; see Corollary 5.4.

This paper is organized as follows. In Section 1, we recall the notion of Hilbert-Kunz multiplicity and its fundamental properties. In Section 2, we prove a generalization of Theorem 1. Section 3 is devoted to study the Hilbert-Kunz multiplicities of the Segre product of polynomial rings. Actually, using the notion of the Stirling numbers of the second kind, we give another proof for the formula in [2] about Hilbert-Kunz multiplicity of Segre products. Also, we show that the constant $c(d)$ appeared in Theorem 2 is given as the Hilbert-Kunz multiplicity of the Rees algebra $R(\mathfrak{m})$ over a polynomial ring with $d$-variables. In Section 4, we prove Theorem 2. Finally, in Section 5, we give an explicit formula for the Hilbert-Kunz multiplicities of the Rees algebra $R(\mathfrak{m})$ over the Veronese subring $k[x_1, \ldots, x_d]^{(c)}$.

**Acknowledgement.** The authors are grateful to Professor Kei-ichi Watanabe for giving them valuable informations. The authors would also like to appreciate the referee's reading this manuscript patiently.

## 1. PRELIMINARIES

In this section, we recall several definitions and fundamental properties about Hilbert-Kunz multiplicity which are needed later; see also e.g. [1,2,4,9,12,24].

Throughout this paper, we use the following notation: For a finitely generated $A$-module $M$, $l_A(M)$ (resp. $\mu_A(M)$) denotes the length of $M$ (resp. the minimal number of generator of $M$).

**1.1. Hilbert-Kunz multiplicity.** First, we recall the notion of the Hilbert-Kunz multiplicity. Let $(A, \mathfrak{m}, k)$ be a local ring of characteristic $p > 0$ with $d := \dim A \geq 1$. Let $I \subseteq A$ be an $\mathfrak{m}$-primary ideal and $M$ a finitely generated $A$-module. For each $q = p^e$, we denote by $I^{[q]}$ the ideal generated by the $q$-th powers of the elements of $I$. Then there exists a positive real constant $C$ such that

$$l_A(M/I^{[q]}M) = Cq^d + O(q^{d-1}) \qquad \text{for all large } q = p^e.$$

Thus we can define the *Hilbert-Kunz multiplicity* of $M$ with respect to $I$ as follows:

$$(1.1.1) \qquad \mathrm{e}_{\mathrm{HK}}(I, M) := \lim_{e \to \infty} \frac{l_A(M/I^{[q]}M)}{q^d}.$$

By definition, we have $\mathrm{e}_{\mathrm{HK}}(I) := \mathrm{e}_{\mathrm{HK}}(I, A)$ and $\mathrm{e}_{\mathrm{HK}}(A) := \mathrm{e}_{\mathrm{HK}}(\mathfrak{m})$. Moreover, if $A$ is a graded ring with unique homogeneous maximal ideal $P$, then we define as $\mathrm{e}_{\mathrm{HK}}(I, A) := \mathrm{e}_{\mathrm{HK}}(IA_P, A_P)$. See also [12,14,15,17].

We also recall the definition of usual multiplicity. For an $\mathfrak{m}$-primary ideal $I$ in $A$, we define the multiplicity $e(I, M)$ with respect to $I$ as follows:

$$(1.1.2) \qquad e(I, M) := \lim_{n \to \infty} \frac{d!}{n^d} l_A(M/I^n M);$$

see [16], [18] or [3] for further details. By definition, we have $e(I) := e(I, A)$ and $e(A) := e(\mathfrak{m})$.

**1.2. Rees algebras.** Next, we recall the definition of (extended) Rees algebras. Let $\mathcal{F} := \{F_n\}_{n \in \mathbb{Z}}$ be a filtration of $A$, that is, $\mathcal{F}$ is a set of ideals which satisfy the following conditions:

(a) Each $F_i$ is an ideal of $A$ and $F_i \supseteq F_{i+1}$ for every integer $i$.
(b) $F_i = A$ for all $i \leq 0$ and $\mathfrak{m} \supseteq F_1$.
(c) $F_i F_j \subseteq F_{i+j}$ for all $i$, $j$.

For such a filtration $\mathcal{F}$, one can define the following graded rings:

$$R(\mathcal{F}) := \bigoplus_{n=0}^{\infty} F_n t^n, \quad R'(\mathcal{F}) := \bigoplus_{n \in \mathbb{Z}} F_n t^n$$

and

$$G(\mathcal{F}) := \bigoplus_{n=0}^{\infty} F_n / F_{n+1} = R'(\mathcal{F})/t^{-1} R'(\mathcal{F}).$$

We call $R(\mathcal{F})$ (resp. $R'(\mathcal{F})$, $G(\mathcal{F})$) the *Rees Algebra* (resp. *the extended Rees Algebra, the associated graded ring*) with respect to $\mathcal{F}$.

For an ideal $I$ of $A$, $\{I^n\}_{n \in \mathbb{Z}}$ gives an important example of filtrations of $A$. In this case, we write $R(I)$ (resp. $R'(I), G(I)$) in place of $R(\{I^n\})$ (resp. $R'(\{I^n\}), G(\{I^n\})$) and call it the Rees algebra (resp. the extended Rees algebra, the associated graded ring) of $I$.

In the following, we present several fundamental properties of the Hilbert-Kunz multiplicity. The next lemma gives a relationship between Hilbert-Kunz multiplicities and usual multiplicities.

**Lemma 1.3.** *Let $(A, \mathfrak{m})$ be a local ring, and let $I \subseteq A$ be an $\mathfrak{m}$-primary ideal. Then*

$$\frac{e(I)}{d!} \leq e_{\mathrm{HK}}(I) \leq e(I).$$

*In particular, if $I$ is a parameter ideal, then $e_{\mathrm{HK}}(I) = e(I)$.*

*Remark 1.4.* Recently, Hanes proved that the first inequality is always strict if $\dim A \geq 2$ in [8].

The next result shows that the tight closure $I^*$ of $I$ is the largest ideal containing $I$ having the same Hilbert-Kunz multiplicity as $I$. Recall that an element $x \in A$ is in the tight closure $I^*$ of $I$ if there exists $c \in A \setminus \bigcup_{\mathfrak{p} \in \mathrm{Min}(A)} \mathfrak{p}$ such that $cx^q \in I^{[q]}$ for all sufficiently large $q = p^e$. See [10,12] for details.

**Lemma 1.5.** ([10, Theorem 8.17]) *Let $(A, \mathfrak{m})$ be a local ring, and let $J \subseteq I$ be ideals of $A$. If $I \subseteq J^*$, then $e_{\mathrm{HK}}(J) = e_{\mathrm{HK}}(I)$. If, in addition, $A$ is excellent reduced, equidimensional, then the converse is also true.*

To determine the Hilbert-Kunz multiplicity of a module-finite subring of a regular local ring (or a polynomial ring), then the next lemma will be useful. For instance, applying the next formula to Veronese subring $A = k[(x,y)^e]$ and $B = k[x,y]$, we obtain $e_{\mathrm{HK}}(A) = \frac{e+1}{2}$.

**Lemma 1.6.** ([24, Theorem 2.7]; see also [2]) *Let $(A, \mathfrak{m}) \subset (B, \mathfrak{n})$ be an extension of local domains where $B$ is a finite $A$-module of rank $r$ and $A/\mathfrak{m} = B/\mathfrak{n}$. Then for every $\mathfrak{m}$-primary ideal $I$, we have*

$$\mathrm{e}_{\mathrm{HK}}(I) = \frac{1}{r}\,\mathrm{e}_{\mathrm{HK}}(IB).$$

*In particular, if $B$ is regular, then $\mathrm{e}_{\mathrm{HK}}(I) = \dfrac{1}{r} l(B/IB)$.*

The following lemma enables us to calculate the Hilbert-Kunz multiplicities of binomial hypersurfaces.

**Lemma 1.7** ([4, Theorem 3.1]). *Let $s, t, d_1, \ldots, d_s, e_1, \ldots, e_t$ be positive integers. Put $u = \max\{d_1, \ldots, d_s, e_1, \ldots, e_t\}$. Then the Hilbert-Kunz multiplicity of a binomial hypersurface $A = k[[x_1, \ldots, x_s, y_1, \ldots, y_t]]/(x_1^{d_1} \cdots x_s^{d_s} - y_1^{e_1} \cdots y_t^{e_t})$ is given as follows:*

$$\mathrm{e}_{\mathrm{HK}}(A) = \sum_{j=1}^{s} \sum_{\ell=1}^{t} (-1)^{j+\ell} s_j t_\ell \frac{j\ell}{(j+\ell-1)u^{j+\ell-1}}.$$

*where $s_j$ and $t_l$ are the $j$-th and $l$-th elementary symmetric polynomials in $d_i$'s and $e_i$'s, respectively.*

*Remark 1.8.* The assumption that "$F$ is homogeneous" in [Co, Theorem 3.1] is superfluous.

## 2. PROOF OF THEOREM 1

In this section, we will prove Theorem 2.1 which is a slight generalization of Theorem 1. This theorem enables us to estimate the Hilbert-Kunz multiplicity of the extended Rees algebra with respect to any filtration $\mathcal{F}$ in terms of some invariants (e.g. $\mathrm{e}_{\mathrm{HK}}(A)$ and $e(A)$) of its base ring.

**Theorem 2.1.** *Let $(A, \mathfrak{m})$ be a local ring of characteristic $p > 0$ with $d = \dim A \geq 1$. Let $\mathcal{F} = \{F_n\}$ be a filtration of $A$ such that $R'(\mathcal{F})$ is a Noetherian ring with $\dim R'(\mathcal{F}) = d+1$. Also, let $I$ be an $\mathfrak{m}$-primary ideal such that $I \supseteq F_1$. If we put $\mathfrak{N} = (t^{-1}, I, R'(\mathcal{F})_+)$ and $\mathfrak{M} = (t^{-1}, \mathfrak{m}, R'(\mathcal{F})_+)$, then we have*
  (1) $\mathrm{e}_{\mathrm{HK}}(I) \leq \mathrm{e}_{\mathrm{HK}}(\mathfrak{N}, R'(\mathcal{F}))$.
  (2) $\mathrm{e}_{\mathrm{HK}}(\mathfrak{N}, R'(\mathcal{F})) \leq \mathrm{e}_{\mathrm{HK}}(G(\mathcal{F}))$, *provided that $F_1$ is an $\mathfrak{m}$-primary ideal.*

*Proof.* In the proof, we put $R' := R'(\mathcal{F})$ and $G := G(\mathcal{F})$ for simplicity. Note that $G$ is a Noetherian ring with $\dim G = d$ by the assumption. Also, $[L]_r$ denotes the homogeneous part with degree $r$ for any graded $R(\mathcal{F})$-module $L$.

(1) Fix $q = p^e$. Set $\mathfrak{N}_i = t^{-i}R' + (IR' + R'_+)^{[q]}$ for all $i = 0, 1, \cdots, q$. Considering a filtration $R' = \mathfrak{N}_0 \supseteq \mathfrak{N}_1 \supseteq \cdots \supseteq \mathfrak{N}_q = \mathfrak{N}^{[q]}$, we have

$$\begin{aligned}\mathfrak{N}_{i-1}/\mathfrak{N}_i &\cong \frac{t^{-i+1}R'}{t^{-i+1}R' \cap (t^{-i}R' + (IR' + R'_+)^{[q]})} \\ &\cong \frac{R'}{t^{-1}R' + (IR' + R'_+)^{[q]} : t^{-i+1}}.\end{aligned}$$

Since $(IR' + R'_+)^{[q]} : t^{-i+1} \subseteq (IR' + R'_+)^{[q]} : t^{-q+1}$ for all $i \leq q$, we get

$$l_{R'}(R'/\mathfrak{N}^{[q]}) = \sum_{i=1}^{q} l_{R'}(\mathfrak{N}_{i-1}/\mathfrak{N}_i)$$

**(2.1.1)**
$$\geq q \cdot l_{R'}(R'/t^{-1}R' + (IR' + R'_+)^{[q]} : t^{-q+1}).$$

We now prove the following claim.

**Claim.** $\left[t^{-1}R' + (IR' + R'_+)^{[q]} : t^{-q+1}\right]_r \subseteq (F_{r+1} + F_r \cap I^{[q]})t^r$ for all integers $r \geq 1$.

Let $at^r$ be an element of $\left[(IR' + R'_+)^{[q]} :_{R'} t^{-q+1}\right]_r$. Then $a \in F_r$ and

$$at^{r-q+1} \in I^{[q]}R' + (R'_+)^{[q]} = \left(I^{[q]} + \sum_{i=1}^{\infty} F_i^{[q]} t^{iq}\right) \sum_{n \in \mathbb{Z}} F_n t^n.$$

Hence one can get

$$t^{-r}\left[t^{-1}R' + (IR' + R'_+)^{[q]} : t^{-q+1}\right]_r$$

**(2.1.2)**
$$= F_{r+1} + F_r \cap \left(I^{[q]}F_{r-q+1} + \sum_{i=1}^{\infty} F_i^{[q]} F_{r-(i+1)q+1}\right).$$

By the assumption that $F_i^{[q]} \subseteq F_1^{[q]} \subseteq I^{[q]}$, we have

$$F_r \cap \left(I^{[q]}F_{r-q+1} + \sum_{i=1}^{\infty} F_i^{[q]} F_{r-(i+1)q+1}\right) \subseteq F_r \cap I^{[q]}.$$

This completes the proof of the above claim.

By virtue of Eq.(2.1.1) and Eq.(2.1.2), we have

$$l_{R'}(R'/\mathfrak{N}^{[q]}) \geq q \cdot \sum_{r=0}^{\infty} l_A(F_r/F_{r+1} + I^{[q]} \cap F_r)$$

$$= q \cdot \sum_{r=0}^{\infty} l_A\left(\frac{F_r + I^{[q]}}{F_{r+1} + I^{[q]}}\right) =: (*).$$

By the assumption that $R'(\mathcal{F})$ is Noetherian, we can find an integer $r = r(q)$ such that $F_r \subseteq I^{[q]}$, and thus $(*) = q \cdot l_A(A/I^{[q]})$. It follows that the following inequality holds:

$$\frac{l_{R'}(R'/\mathfrak{N}^{[q]})}{q^{d+1}} \geq \frac{l_A(A/I^{[q]})}{q^d}$$

for all $q = p^e$. Let $e$ tend to $\infty$, and we obtain $e_{HK}(\mathfrak{N}, R') \geq e_{HK}(I)$, as required.

(2) Since $F_1$ is an $\mathfrak{m}$-primary ideal, $(\mathfrak{m}/F_1)G$ is a nilpotent ideal of $G$. Thus $G_+ \subseteq \mathfrak{N}G \subseteq (G_+)^* = \mathfrak{M}G$. Moreover, as $t^{-1} \in \mathfrak{N}$ is a homogeneous non-zero-divisor in $R'$, by virtue of Lemma 1.5 and [24, Proposition 2.13], we get $e_{HK}(G) = e_{HK}(G_+G, G) = e_{HK}(\mathfrak{N}G, G) \geq e_{HK}(\mathfrak{N}R', R')$, as required. $\square$

Note that Theorem 1 easily follows from the following corollary.

**Corollary 2.2.** *Under the same notation as in Theorem* 2.1, *we also assume that $F_1$ is an $\mathfrak{m}$-primary ideal. Then*
  (1) $\mathrm{e}_{\mathrm{HK}}(F_1) \leq \mathrm{e}_{\mathrm{HK}}((t^{-1}, R'(\mathcal{F})_+)R'(\mathcal{F})) \leq \mathrm{e}_{\mathrm{HK}}(G(\mathcal{F}))$.
  (2) $\mathrm{e}_{\mathrm{HK}}(A) \leq \mathrm{e}_{\mathrm{HK}}(R'(\mathcal{F})) \leq \mathrm{e}_{\mathrm{HK}}(G(\mathcal{F}))$.

**Example 2.3.**
  (1) *Let $A = k[A_1] = \oplus_{n \geq 0} A_n$ be a homogeneous $k$-algebra over a field $k$ of characteristic $p > 0$. Put $\mathfrak{m} = A_+$. Then $\mathrm{e}_{\mathrm{HK}}(A) = \mathrm{e}_{\mathrm{HK}}(R'(\mathfrak{m})) = \mathrm{e}_{\mathrm{HK}}(G(\mathfrak{m}))$.*
  (2) *Let $I$ be an ideal of a local ring $(A, \mathfrak{m})$ of characteristic $p > 0$. If $\mathrm{e}_{\mathrm{HK}}(I) = e(I)$, then $\mathrm{e}_{\mathrm{HK}}(G(I)) = e(G(I)) = e(I)$. In particular, if $I$ is a parameter ideal of $A$, then $\mathrm{e}_{\mathrm{HK}}(G(I)) = e(I)$.*

Note that $\mathrm{e}_{\mathrm{HK}}(I) \leq \mathrm{e}_{\mathrm{HK}}(R'(I))$ does not hold in general even if $I$ is a parameter ideal of a Cohen–Macaulay local ring $A$; see Example 5.6.

In Corollary 2.2, even if $\mathcal{F} = \{\mathfrak{m}^n\}$, equalities do not hold in general. Actually, we give an example below which makes equality fail.

**Example 2.4.** (cf. [24, Sect.5]) Let $k$ be any field of characteristic $p > 0$, and let $n$ be any positive integer with $n \geq 2$. Let $A = k[x, y, z]/(xy - z^n)$. Then
$$R'(\mathfrak{m}) \cong k[x, y, z, w]/(xy - z^n w^{n-2})$$
and we have
$$\mathrm{e}_{\mathrm{HK}}(A) = 2 - \frac{1}{n}, \qquad \mathrm{e}_{\mathrm{HK}}(R'(\mathfrak{m})) = 2 - \frac{2(n+1)}{3n^2}.$$
Moreover, when $n \geq 3$, we have $G(\mathfrak{m}) \cong k[x, y]/(xy)$; thus $\mathrm{e}_{\mathrm{HK}}(G(\mathfrak{m})) = 2$.

*Proof.* Because $A$ and $R'$ are both hypersurfaces, we can apply Lemma 1.7 (the method of Conca) to them. □

*Discussion 2.5.* Let $A = k[[x, y, z]]/(f)$ be a two-dimensional hypersurface of multiplicity 2. If $A$ is *not* F-rational, then $\mathrm{e}_{\mathrm{HK}}(A) = e(A)$. Hence $\mathrm{e}_{\mathrm{HK}}(R'(\mathfrak{m})) = 2$ by Theorem 2.1. Therefore if one wants to get $\mathrm{e}_{\mathrm{HK}}(R'(\mathfrak{m}))$ completely, one may assume that $A$ is F-rational.

Further, assume that $k$ is algebraically closed field. Then the $\mathfrak{m}$-adic completion of any F-rational hypersurface is isomorphic to a local ring which is called "F-rational double point" defined by either one of the following equations:

$$\begin{aligned}
A_n &: xy + z^{n+1} & (n \geq 1), \quad &\text{where } p \geq 2 \\
D_n &: x^2 + yz^2 + y^{n-1} & (n \geq 4), \quad &\text{where } p \geq 3 \\
E_6 &: x^2 + y^3 + z^4, & &\text{where } p \geq 5 \\
E_7 &: x^2 + y^3 + yz^3, & &\text{where } p \geq 5 \\
E_8 &: x^2 + y^3 + z^5, & &\text{where } p \geq 7
\end{aligned}$$

By virtue of [13, Corollary 4.4], $R'(\mathfrak{m})$ is also an F-rational hypersurface for any local ring in the above list. This implies that $\mathrm{e}_{\mathrm{HK}}(R'(\mathfrak{m})) < 2$. But it seems to be difficult to determine $\mathrm{e}_{\mathrm{HK}}(R'(\mathfrak{m}))$ except $(A_n)$. □

We also pose the following question.

**Question 2.6.** Let $A = k[[x, y, z]]/(xy - z^{n+1}) \cong k[[s^{n+1}, st, t^{n+1}]]$, the rational double point of type $(A_n)$. How about $e_{HK}(R(\mathfrak{m}))$?

## 3. CALCULUS OF THE HILBERT-KUNZ MULTIPLICITY OF THE SEGRE PRODUCT

In this section, we will show that the constant $c(d)$ appeared in Theorem 2 is equal to the Hilbert-Kunz multiplicity of the Rees algebra with respect to the maximal ideal over a polynomial ring (or a regular local ring) of dimension $d$. Namely, we have

$$(3.1) \qquad e_{HK}(R(\mathfrak{m})) = c(d) = d\left(\frac{1}{2} + \frac{1}{(d+1)!}\right),$$

where $A = k[x_1, \ldots, x_d]$ and $\mathfrak{m} = (x_1, \ldots, x_d)A$. Then since the Rees algebra $R(\mathfrak{m})$ is isomorphic to the Segre product $S_{2,d} := k[x_1, x_2] \# k[y_1, \ldots, y_d]$, it is enough to calculate the Hilbert-Kunz multiplicity of the Segre product of polynomial rings in general. In fact, Buchweitz et.al showed the following formula in [2]:

$$(3.2) \quad e_{HK}(S_{c,d}) \bigg/ \binom{c+d}{c}$$
$$= \frac{(c+1)^{c+d+1}}{(c+d+1)!} - \sum_{0 \leq i, k < j \leq c} \frac{(-1)^{i+k}}{(c+d)!} \binom{d+1}{j-i} \binom{c+1}{j-k} \int_0^1 (u+i)^d (u+k)^c du.$$

But we hope that this formula will become more clear! In fact, we prove the following theorem in terms of "the Stirling number of the second kind".

**Theorem 3.3.** (cf. [2],[5]) Let $R = k[x_1, \ldots, x_c] \# k[y_1, \ldots, y_d]$ be the Segre product of polynomial rings over a field $k$ with $c, d$ variables, respectively, where $c \leq d$. Then

$$e_{HK}(R) = \frac{d!}{(c+d-1)!} S(c+d-1, d)$$
$$- \frac{1}{(c+d-1)!} \sum_{0 < j < i \leq c} \binom{c}{i}\binom{d}{j}(-1)^{c-i+j}(i-j)^{c+d-1},$$

where $S(n, k)$ denotes the Stirling number of the second kind; see below.

In particular, if $c = 2$,
$$e_{HK}(R) = d\left(\frac{1}{2} + \frac{1}{(d+1)!}\right)$$

For example,
$$e_{HK}(S_{2,2}) = \frac{4}{3}, \quad e_{HK}(S_{3,3}) = \frac{39}{20}, \quad e_{HK}(S_{4,4}) = \frac{899}{315}, \quad e_{HK}(S_{5,5}) = \frac{151205}{36288},$$
$$e_{HK}(S_{6,6}) = \frac{10114043}{1663200}, \quad e_{HK}(S_{2,3}) = \frac{13}{8}, \quad e_{HK}(S_{3,4}) = \frac{889}{360}.$$

We note that $e_{HK}(S_{4,4}) = \frac{899}{315}$ is different from the value in [BCP,2.2.3].

Before proving this theorem, we recall the notion of Stirling numbers and gather several properties which are needed later. See e.g. [21, Chapter 1, §1.4] for details.

**Definition 3.4.** For a given natural number $k$ and a polynomial $f(x) \in \mathbb{C}[x]$, we denote $f(x)f(x-1)\cdots f(x-k+1)$ by $f^{\underline{k}}$. And we put $f^{\underline{0}} = 1$. Then since $\{x^{\underline{k}}\}_{k \geq 0}$ forms a basis of the vector space $\mathbb{C}[x]$ over $\mathbb{C}$, there uniquely exist integers $S(n, k)$ where $n \geq 0$ such that

$$x^n = \sum_{k=0}^{\infty} S(n,k) x^{\underline{k}} = \sum_{k=0}^{\infty} S(n,k)\, x(x-1)\cdots(x-k+1).$$

Then $S(n, k)$ is called the *Stirling number of the second kind*.

One can easily get $S(n, k) = 0$ if $k > n$ and $S(0, 0) = 1$ by definition. Also, for all $n \geq 1$, one has

$$S(n,0) = 0, \quad S(n,1) = 1, \quad S(n,2) = 2^{n-1} - 1, \quad S(n,n) = 1, \quad S(n, n-1) = \binom{n}{2}.$$

**Fact 3.5.** (cf. [21]) *The Stirling number of the second kind admits the following characterizations*:

(1) $S(n,k)$ is equal to the number of partitions of the set $[n] := \{1, \ldots, n\}$ into $k$ blocks.
(2) The Stirling numbers of the second kind satisfy the following recurrence:

$$S(n,k) = k \cdot S(n-1,k) + S(n-1, k-1), \quad S(0,0) = 1, \quad S(n,0) = S(0,k) = 0.$$

(3) $S(n,k)$ admits the following exponential generating function:

$$\sum_{n \geq k} S(n,k) \frac{x^n}{n!} = \frac{1}{k!}(e^x - 1)^k.$$

*In particular,*

$$S(n,k) = \frac{1}{k!} \sum_{i=0}^{k} (-1)^{k-i} \binom{k}{i} i^n.$$

Table 1. Stirling numbers of the second kind S(n,k).

| $n \setminus k$ | 1 | 2 | 3 | 4 | 5 | 6 | 7 | 8 | 9 | 10 |
|---|---|---|---|---|---|---|---|---|---|---|
| 1 | 1 | 0 | 0 | 0 | 0 | 0 | 0 | 0 | 0 | 0 |
| 2 | 1 | 1 | 0 | 0 | 0 | 0 | 0 | 0 | 0 | 0 |
| 3 | 1 | 3 | 1 | 0 | 0 | 0 | 0 | 0 | 0 | 0 |
| 4 | 1 | 7 | 6 | 1 | 0 | 0 | 0 | 0 | 0 | 0 |
| 5 | 1 | 15 | 25 | 10 | 1 | 0 | 0 | 0 | 0 | 0 |
| 6 | 1 | 31 | 90 | 65 | 15 | 1 | 0 | 0 | 0 | 0 |
| 7 | 1 | 63 | 301 | 350 | 140 | 21 | 1 | 0 | 0 | 0 |
| 8 | 1 | 127 | 966 | 1701 | 1050 | 266 | 28 | 1 | 0 | 0 |
| 9 | 1 | 255 | 3025 | 7770 | 6951 | 2646 | 462 | 36 | 1 | 0 |
| 10 | 1 | 511 | 9330 | 34105 | 42525 | 22827 | 5880 | 750 | 45 | 1 |

In the following, let $A_d := k[x_1, \ldots, x_d]$ be a polynomial ring over a field $k$ with $d$ variables and put $\mathfrak{m} = (x_1, \ldots, x_d)A$. Also, if one set

$$(3.6.1) \quad \alpha_{d,n} := l_A(\mathfrak{m}^n/\mathfrak{m}^{n+1}) = \binom{n+d-1}{d-1}, \quad \alpha_{d,n,q} := l_A(\mathfrak{m}^n/\mathfrak{m}^{n-q}\mathfrak{m}^{[q]} + \mathfrak{m}^{n+1})$$

for all integers $q = p^e$ and $m$, then one can easily obtain that

$$(3.6.2) \quad \alpha_{d,n,q} = \sum_{i=0}^{d}(-1)^i \binom{d}{i} \alpha_{d,n-iq}.$$

Actually, $\alpha_{d,n,q}$ is the number of monomials of degree $n$ which appears in the polynomial $\prod_{i=1}^{d}(1 + x_i + x_i^2 + \cdots + x_i^{q-1})$.

From now on, we will prove Theorem 3.3.

**Lemma 3.7.** (See [7], [2,5]) *Let $R = A_c \# A_d$ be the Segre product with $2 \leq c \leq d$. Then we have*

(1) $\dim R = c + d - 1$.
(2) *The Hilbert-Kunz multiplicity of $R$ can be calculated by the following formula:*

$$\mathrm{e}_{\mathrm{HK}}(R) = \lim_{q \to \infty} \frac{1}{q^{c+d-1}} \sum_{n=0}^{d(q-1)} \alpha_{c,n}\alpha_{d,n,q}$$

$$(3.7.1) \quad + \lim_{q \to \infty} \frac{1}{q^{c+d-1}} \sum_{n=0}^{c(q-1)} \alpha_{c,n,q}\alpha_{d,n} - \lim_{q \to \infty} \frac{1}{q^{c+d-1}} \sum_{n=0}^{c(q-1)} \alpha_{c,n,q}\alpha_{d,n,q}.$$

To calculate the first term and the second term in Eq.(3.7.1), we need the following lemma. Notice that the main calculation of the following two lemmata follows from Lemma 3.10 below.

**Lemma 3.8.** *Under the above notation, for all positive integers $c$, $d$, we have*

$$\lim_{q \to \infty} \frac{1}{q^{c+d-1}} \sum_{n=0}^{c(q-1)} \alpha_{c,n,q}\alpha_{d,n} = \frac{c!}{(c+d-1)!}S(c+d-1, c).$$

*Proof.* It is known that $\alpha_{d,n} = \frac{1}{(d-1)!}n^{d-1} + $(lower term) if $n \gg 0$. Hence we may assume that $\alpha_{d,n}$ is a polynomial of degree $d-1$ with the leading coefficient $\frac{1}{(d-1)!}$. Then by Eq.(3.6.1), (3.6.2), we have

$$\lim_{q \to \infty} \frac{1}{q^{c+d-1}} \sum_{n=0}^{c(q-1)} \alpha_{c,n,q}\alpha_{d,n}$$

$$(3.8.1) \quad = \lim_{q \to \infty} \frac{1}{q^{c+d-1}} \sum_{i=0}^{c} \binom{c}{i}(-1)^i \sum_{n=iq}^{c(q-1)} \frac{(n-iq)^{c-1}}{(c-1)!} \cdot \frac{n^{d-1}}{(d-1)!}$$

$$= \frac{1}{(c-1)!(d-1)!} \sum_{i=0}^{c} \binom{c}{i}(-1)^i \int_{i}^{c} (t-i)^{c-1}t^{d-1}\, dt.$$

Applying Lemma 3.10 as $j = 0$, we have

$$\text{LHS of (3.8.1)} = \frac{1}{(c+d-1)!} \sum_{i=0}^{c} \binom{c}{i} (-1)^{c-i} i^{c+d-1}$$

$$= \frac{c!}{(c+d-1)!} S(c+d-1, c),$$

where the second equality follows from Fact 3.5 (3). □

Furthermore, we need the following lemma to calculate the last term in Eq.(3.7.1).

**Lemma 3.9.** *Under the same notation as in the previous lemma, for integers $0 < c \leq d$, we have*

$$(3.9.1) \quad \lim_{q \to \infty} \frac{1}{q^{c+d-1}} \sum_{n=0}^{c(q-1)} \alpha_{c,n,q} \alpha_{d,n,q} = \frac{c!}{(c+d-1)!} S(c+d-1, c)$$

$$+ \frac{1}{(c+d-1)!} \sum_{0 < j < i \leq c} \binom{c}{i}\binom{d}{j} (-1)^{c-i+j} (i-j)^{c+d-1}.$$

*Proof.* By the similar argument as in the proof of Lemma 3.8, we have

$$\lim_{q \to \infty} \frac{1}{q^{c+d-1}} \sum_{n=0}^{c(q-1)} \alpha_{c,n,q} \alpha_{d,n,q}$$

$$= \frac{1}{(c-1)!(d-1)!} \sum_{i=0}^{c} \sum_{j=0}^{c} \binom{c}{i}\binom{d}{j} (-1)^{i+j} \int_{\max\{i,j\}}^{c} (t-i)^{c-1}(t-j)^{d-1} \, dt$$

$$= \frac{1}{(c-1)!(d-1)!} \left[ \sum_{i=0}^{c} \sum_{j=0}^{c} \int_{j}^{c} - \sum_{0 \leq j < i \leq c} \int_{j}^{i} \right] \binom{c}{i}\binom{d}{j} (-1)^{i+j} (t-i)^{c-1}(t-j)^{d-1} \, dt.$$

By Lemma 3.10, we have

$$\text{LHS of (3.9.1)} = \frac{1}{(c-1)!(d-1)!} \sum_{0 \leq j < i \leq c} \binom{c}{i}\binom{d}{j} (-1)^{i+j+c} (i-j)^{c+d-1} B(c,d)$$

$$= \frac{1}{(c+d-1)!} \left( \sum_{\substack{j=0, \\ 1 \leq i \leq c}} + \sum_{0 < j < i \leq c} \right) \binom{c}{i}\binom{d}{j} (-1)^{c-i+j} (i-j)^{c+d-1}.$$

The required equality easily follows from Fact 3.5(3) and the above equality. □

*Proof of Theorem 3.3.* By the above lemmata, we have

$$\mathrm{e}_{\mathrm{HK}}(R) = \frac{c!}{(c+d-1)!}S(c+d-1,c) + \frac{d!}{(c+d-1)!}S(c+d-1,d)$$
$$- \frac{c!}{(c+d-1)!}S(c+d-1,c)$$
$$- \frac{1}{(c+d-1)!} \sum_{0<j<i\leq c} \binom{c}{i}\binom{d}{j}(-1)^{c-i+j}(i-j)^{c+d-1}.$$

This yields the required equality. □

**Lemma 3.10.** *Let $i$, $j$, and $c$ be integers. Then*

(1) $\displaystyle\sum_{i=0}^{c} \binom{c}{i}(-1)^i \int_j^c (t-i)^{c-1}(t-j)^{d-1}\, dt = 0.$

(2) $\displaystyle\int_j^i (t-i)^{c-1}(t-j)^{d-1}\, dt = (-1)^{c-1}(i-j)^{c+d-1}B(c,d),$

*where*

$$B(c,d) = \frac{(c-1)!(d-1)!}{(c+d-1)!}.$$

*Proof.* (1) Since $\displaystyle\sum_{i=0}^{c}(-1)^i\binom{c}{i}i^n = 0$ for all $0 \leq n \leq c-1$ by Fact 3.5(3),

$$\text{LHS of (1)} = \int_j^c (t-j)^{d-1}\left\{\sum_{i=0}^{c}\binom{c}{i}(-1)^i(t-i)^{c-1}\right\} dt = 0.$$

(2) Putting $x = \frac{t-j}{i-j}$, we have $1-x = \frac{i-t}{i-j}$ and $dx = \frac{1}{i-j}\, dt$. Then

$$\text{LHS of (2)} = \int_0^1 (-1)^{c-1}(i-j)^{c+d-1}(1-x)^{c-1}x^{d-1}\, dx$$
$$= (-1)^{c-1}(i-j)^{c+d-1}\int_0^1 (1-x)^{c-1}x^{d-1}\, dx.$$

Also, since it is known that

$$\int_0^1 (1-x)^{c-1}x^{d-1}\, dx = B(c,d) = \frac{(c-1)!(d-1)!}{(c+d-1)!}$$

by the property of the beta function, we get the required equality. □

# 4. PROOF OF THEOREM 2

In this section, we will prove Theorem 2. First we prove the right-hand side inequality in Theorem 2.

**Theorem 4.1.** *Let $(A, \mathfrak{m})$ be a local ring of characteristic $p > 0$ with $d = \dim A \geq 1$. Then for any $\mathfrak{m}$-primary ideal $I$, we have*

$$\mathrm{e}_{\mathrm{HK}}(R(I)) \leq c(d) \cdot e(I),$$

*where $c(d) = d(\frac{1}{2} + \frac{1}{(d+1)!})$.*

*Moreover, equality holds if and only if $\mathrm{e}_{\mathrm{HK}}(A) = e(I)$. When this is the case, $\mathrm{e}_{\mathrm{HK}}(A) = e(A)$ and $\mathrm{e}_{\mathrm{HK}}(I) = e(I)$.*

We now begin our proof of the above theorem by giving the next well-known lemma.

**Lemma 4.2.** *Let $(A, \mathfrak{m})$ be any local ring. Let $I$ be an $\mathfrak{m}$-primary ideal of $A$, and put $R := R(I) = A[It]$. For any prime $\mathfrak{p}$ in $A$, we set $\mathfrak{p}^* = \mathfrak{p}A[t] \cap R(I)$. Put $\mathrm{Assh}(A) := \{\mathfrak{p} \in \mathrm{Spec}\, A \mid \dim A/\mathfrak{p} = \dim A\}$. Then the following statements hold.*

(1) $\mathrm{Assh}(R) = \{\mathfrak{p}^* \in \mathrm{Spec}\, R \mid \mathfrak{p} \in \mathrm{Assh}(A)\}$.
(2) $R/\mathfrak{p}^* \cong R(I + \mathfrak{p}/\mathfrak{p}, A/\mathfrak{p})$, *the Rees algebra of the ideal $I + \mathfrak{p}/\mathfrak{p}$ in $A/\mathfrak{p}$*.
(3) $R_{\mathfrak{p}^*} \cong A_{\mathfrak{p}}(t)$ *if $\mathfrak{p} \not\supseteq I$. In particular, if $\mathfrak{p}$ is a minimal prime ideal in $A$, then $l_{A_{\mathfrak{p}}}(A_{\mathfrak{p}}) = l_{R_{\mathfrak{p}}^*}(R_{\mathfrak{p}^*})$.*

By the above lemma and [24,(2.3)], we also obtain the next lemma.

**Lemma 4.3.** *Under the same notation as in the previous lemma, we further assume that $A$ has characteristic $p > 0$ and let $L \subseteq A$ be an ideal with $I \subseteq L \subseteq \mathfrak{m}$. Then we get the following formula.*

$$\mathrm{e}_{\mathrm{HK}}((L, It)R(I)) = \sum_{\mathfrak{p} \in \mathrm{Assh}(A)} \mathrm{e}_{\mathrm{HK}}((\overline{L}, \overline{I}t)\overline{A}[\overline{I}t]) \cdot l_{A_{\mathfrak{p}}}(A_{\mathfrak{p}}),$$

*where $\overline{A} = A/\mathfrak{p}$, $\overline{I} = I + \mathfrak{p}/\mathfrak{p}$ and $\overline{L} = L + \mathfrak{p}/\mathfrak{p}$.*

*Proof.* Applying [24,(2.3)] to $R := R(I)$, we get

$$\begin{aligned}
\mathrm{e}_{\mathrm{HK}}((L, It)R) &= \sum_{P \in \mathrm{Assh}(R)} \mathrm{e}_{\mathrm{HK}}((L, It)R/P) \cdot l_{R_P}(R_P) \\
&= \sum_{\mathfrak{p} \in \mathrm{Assh}(A)} \mathrm{e}_{\mathrm{HK}}((L, It)R/\mathfrak{p}^*) \cdot l_{R_{\mathfrak{p}^*}}(R_{\mathfrak{p}^*}) \\
&= \sum_{\mathfrak{p} \in \mathrm{Assh}(A)} \mathrm{e}_{\mathrm{HK}}((\overline{L}, \overline{I}t)\overline{A}[\overline{I}t]) \cdot l_{A_{\mathfrak{p}}}(A_{\mathfrak{p}}),
\end{aligned}$$

as required. □

Using these lemmata, we prove the following proposition.

**Proposition 4.4.** *Let $(A, \mathfrak{m})$ be a local ring of characteristic $p > 0$. Suppose that $I, J$ are $\mathfrak{m}$-primary ideals such that $J$ is a reduction of $I$, that is, $J \subseteq I$ and $I^{n+1} = JI^n$ for some non-negative integer $n$. Then*

$$e_{\mathrm{HK}}(R(I)) \leq e_{\mathrm{HK}}(R(J)).$$

*Proof.* By virtue of the previous lemma, we may assume that $A$ is a local domain. Let $\mathfrak{M}$ (resp. $\mathfrak{N}$) denote the homogeneous maximal ideal of $R(I)$ (resp. $R(J)$). Then $R(I)_{\mathfrak{M}}$ is a local domain which is module-finite over $R(J)_{\mathfrak{N}}$. Further, as $R(I)_{\mathfrak{M}}$ and $R(J)_{\mathfrak{N}}$ have the same fraction field, by virtue of Lemma 1.6, we have $e_{\mathrm{HK}}(R(J)) = e_{\mathrm{HK}}(\mathfrak{N}R(I), R(I)) \geq e_{\mathrm{HK}}(R(I))$. □

As a corollary of the above proposition, we get the following.

**Corollary 4.5.** *Under the same notation as in* Proposition 4.4, *the following statements hold.*
  (1) $e_{\mathrm{HK}}((I, It)R(I)) \leq e_{\mathrm{HK}}((J, Jt)R(J))$.
  (2) *In (1), equality holds if and only if $e_{\mathrm{HK}}(I) = e_{\mathrm{HK}}(J)$.*

*Proof.* One can prove (1) by the similar argument as in Lemma 4.3 and Proposition 4.4. To see (2), we may assume that $A$ is a complete local domain. Put $R := R(I)$. By Lemma 1.5 and Lemma 1.6, equality holds in (1) if and only if $(I, It)R \subseteq ((J, Jt)R)^*$. On the other hand, we know that $e_{\mathrm{HK}}(I) = e_{\mathrm{HK}}(J)$ if and only if $I \subseteq J^*$. Hence it suffices to show that $(I, It)R \subseteq ((J, Jt)R)^*$ if and only if $I \subseteq J^*$.

First, suppose that $(I, It)R \subseteq ((J, Jt)R)^*$. For any $a \in I$, if we regard $a$ as an element of $R$, then $a \in ((J, Jt)R)^*$. Thus we can take $c = c_r t^r + c_{r+1} t^{r+1} + \cdots + c_s t^s \in R$ ($c_r \neq 0$) such that $ca^q \in (J^{[q]}, J^{[q]}t^q)R$ for all $q = p^e$. In particular, we have $c_r a^q \in J^{[q]}I^r + J^{[q]}I^{q-r} \subseteq J^{[q]}$ for all $q = p^e$; hence $a \in J^*$ as required.

Next, suppose $I \subseteq J^*$. To see $(I, It)R \subseteq ((J, Jt)R)^*$, it is enough to prove $It \subseteq ((J, Jt)R)^*$. For any $a \in I$, we can take a non-zero element $c$ such that $ca^q \in J^{[q]}$ for all $q = p^e$. Since $c(at)^q = (ca^q)t^q \in J^{[q]}t^q \subseteq ((J, Jt)R)^{[q]}$ for all $q = p^e$, we get $at \in ((J, Jt)R)^*$. This completes the proof of the corollary. □

*Proof of Theorem* 4.1. First we prove the inequality stated as above. By Lemma 4.3 and Proposition 4.4, we may assume that $A$ is a complete local domain and $I$ is a parameter ideal of $A$. Say $I = (a_1, \ldots, a_d)A$.

Take a coefficient field $k \subseteq A$ and set $B := k[[a_1, \ldots, a_d]] \subseteq A$. Then $B$ is a complete regular local ring with $\dim B = d$ and $A$ is a finitely generated $B$-module. Also, $[Q(A) : Q(B)] = e(I)$, where $Q(A)$ (resp. $Q(B)$) denotes the fraction field of $A$ (resp. $B$).

Now consider the Rees algebra $R := R(I) = A[It]$. Put $\mathfrak{m}_B = (a_1, \ldots, a_d)B$. Then $R(I)$ is a finitely generated $R(\mathfrak{m}_B)$-module and

$$[Q(R(I)) : Q(R(\mathfrak{m}_B))] = [Q(A)(t) : Q(B)(t)] = [Q(A) : Q(B)] = e(I).$$

By Lemma 1.6 and Theorem 3.3, we have

$$e_{\mathrm{HK}}((I, It)R(I)) = [Q(R(I)) : Q(R(\mathfrak{m}_B))] \cdot e_{\mathrm{HK}}((\mathfrak{m}_B, \mathfrak{m}_B t)R(\mathfrak{m}_B)) = e(I) \cdot c(d).$$

Thus $e_{HK}(R(I)) \leq e_{HK}((I, It)R(I)) = c(d) \cdot e(I)$, as required.

Next, suppose that $e_{HK}(R(I)) = c(d) \cdot e(I)$ holds. To see $e_{HK}(A) = e(I)$, we may assume that $A$ is a complete local domain. Also, let $J$ be a minimal reduction of $I$. Then since $e_{HK}(R(J)) \geq e_{HK}(R(I)) = c(d) \cdot e(J) = e_{HK}((J, Jt)R(J))$ by assumption, we get $(\mathfrak{m}, Jt) \subseteq (J, Jt)^*$ in $R(J)$ by Lemma 1.5. By the similar argument as in the proof of Corollary 4.5, we get $\mathfrak{m} \subseteq J^*$. Hence we have $e_{HK}(A) = e_{HK}(J) = e(J) = e(I)$ as required.

Conversely, suppose that $e_{HK}(A) = e(I)$. To see $e_{HK}(R(I)) = c(d) \cdot e(I)$, we may also assume that $A$ is a complete local domain. Then by assumption we obtain $\mathfrak{m} = J^*$. Then $J$ is a reduction of $\mathfrak{m}$ and thus by Corollary 4.5, we get

$$e_{HK}(R(\mathfrak{m})) = e_{HK}((I, It)R(I)) = e_{HK}((J, Jt)R(J)) = c(d) \cdot e(I).$$

Hence $e_{HK}(R(I)) = c(d) \cdot e(I)$, as required. $\square$

*Remark 4.6.* For any $\mathfrak{m}$-primary ideal $I$ of a local ring, we have

$$l_R(R/(\mathfrak{m}, It)^{[q]}) = \sum_{n=0}^{q-1} l_A(I^n/\mathfrak{m}^{[q]}I^n) + \sum_{n=q}^{\infty} l_A(I^n/\mathfrak{m}^{[q]}I^n + I^{[q]}I^{n-q})$$

$$l_R(R/(I, It)^{[q]}) = \sum_{n=0}^{q-1} l_A(I^n/I^{[q]}I^n) + \sum_{n=q}^{\infty} l_A(I^n/I^{[q]}I^{n-q}).$$

for all $q = p^e$.

The following example will be generalized in Section 5.

**Example 4.7 (K.-i. Watanabe).** Let $A = k[[(x, y)^e]]$ and put $\mathfrak{m} = (x, y)^e A$. Then we have

$$e_{HK}(R(\mathfrak{m})) = e + \frac{1}{3e} \leq c(2) \cdot e(A) = \frac{4}{3}e.$$

*Proof.* Put $R := R(\mathfrak{m})$, $\mathfrak{M} := \mathfrak{m}R + R_+ = (\mathfrak{m}, \mathfrak{m}t)$. Then we get

$$l_R(R/\mathfrak{M}^{[q]}) = \sum_{n=0}^{q-1} l_A(\mathfrak{m}^n/\mathfrak{m}^{[q]}\mathfrak{m}^n) + \sum_{n=q}^{2q-1} l_A(\mathfrak{m}^n/\mathfrak{m}^{[q]}\mathfrak{m}^{n-q})$$

$$= 2q \cdot l_A(A/\mathfrak{m}^{[q]}) - \sum_{n=0}^{q-1} l_A(A/\mathfrak{m}^n) - \sum_{n=0}^{q-1} l_A(A/\mathfrak{m}^{n+q})$$

$$+ 2\sum_{n=0}^{q-1} l_A(\mathfrak{m}^{[q]}/\mathfrak{m}^{[q]}\mathfrak{m}^n).$$

Since $l_A(A/\mathfrak{m}^n) = \frac{e}{2}n^2 + \left(1 - \frac{e}{2}\right)n$ for all $n \geq 0$, we get

$$\sum_{n=0}^{q-1} l_A(A/\mathfrak{m}^n) = \frac{e}{6}q^3 + O(q^2), \qquad \sum_{n=0}^{q-1} l_A(A/\mathfrak{m}^{n+q}) = \frac{7}{6}eq^3 + O(q^2).$$

On the other hand, since we have $e_{HK}(A) = \frac{e+1}{2}$, we get

$$2 \cdot l_A(A/\mathfrak{m}^{[q]}) = 2 e_{HK}(A) q^3 + O(q^2) = (e+1) q^3 + O(q^2).$$

Thus it suffices to prove the following claim.

**Claim.** $\sum_{n=0}^{q-1} l_A(\mathfrak{m}^{[q]}/\mathfrak{m}^{[q]}\mathfrak{m}^n) = \left( \frac{2}{3} e - \frac{1}{2} + \frac{1}{6e} \right) q^3 + O(q^2).$

Actually, we get

$$\text{LHS} = \sum_{n=0}^{q-1} \sum_{t=0}^{n-1} \mu_A(\mathfrak{m}^{[q]}\mathfrak{m}^t)$$

$$= \sum_{n=1}^{q} \sum_{t=0}^{n-1} (e+1)(te+1) - \sum_{n=q/e+1}^{n-1} \sum_{t=q/e}^{n-1} e(te - q) + O(q^2)$$

$$= \frac{e(e+1)}{6} q^3 - \sum_{n=q/e+1}^{q} \left[ \frac{e^2}{2} \left\{ n^2 - \left(\frac{q}{e}\right)^2 \right\} - eq \left( n - \frac{q}{e} \right) \right] + O(q^2)$$

$$= \frac{e(e+1)}{6} q^3 - \sum_{n=q/e+1}^{q} \left( \frac{e^2}{2} n^2 - eqn + \frac{q^2}{2} \right) + O(q^2).$$

The assertion follows easily from the above claim. $\square$

In the rest of this section, we will prove the left-hand side inequality in Theorem 2.

**Theorem 4.8.** *Let $(A, \mathfrak{m})$ be a local ring of characteristic $p > 0$ with $d = \dim A \geq 1$. Then*

$$e_{HK}(A) \leq e_{HK}(R'(\mathfrak{m})) \leq e_{HK}(R(\mathfrak{m})).$$

*Remark 4.9.* According to Theorem 2.1, we have $e_{HK}(A) \leq e_{HK}(R'(I))$ for any $\mathfrak{m}$-primary ideal $I$. However, in general, $e_{HK}(R'(I)) \leq e_{HK}(R(I))$ does not necessarily hold; see e.g. Example 5.6.

We begin the proof of the above theorem by giving the next lemma.

**Lemma 4.10.** *Under the same notation as in* Theorem 4.8,

$$\liminf_{q \to \infty} \frac{1}{q^{d+1}} \sum_{n=0}^{q-1} l_A(\mathfrak{m}^n/\mathfrak{m}^{[q]}\mathfrak{m}^n) \geq e_{HK}(A).$$

*Proof.* Put $R = R(\mathfrak{m})$, $G = G(\mathfrak{m})$ and $\mathfrak{M} = (\mathfrak{m}, \mathfrak{m}t)R$. In order to prove this lemma, it suffices to show that

$$\liminf_{q \to \infty} \frac{f(q)}{q^{d+1}} \geq 0, \quad \text{where} \quad f(q) := \sum_{n=0}^{q-1} l_A(\mathfrak{m}^n/\mathfrak{m}^{[q]}\mathfrak{m}^n) - q \cdot l_A(A/\mathfrak{m}^{[q]}).$$

First, we prove the following claim.

**Claim.** $\mu_A(\mathfrak{m}^{[q]}\mathfrak{m}^t) \geq l_A([\overline{G}]_t)$ for every $t \geq 1$, where $\overline{G} = G/H^0_{\mathfrak{M}}(G)$ and $[\overline{G}]_t$ denotes the homogeneous part of $\overline{G}$ with degree $t$.

Take an element $a \in \mathfrak{m} \setminus \mathfrak{m}^2$ such that $\mathrm{in}(a) \in \overline{G}_1$ is a non-zero-divisor. Fix $t \in \mathbb{N}$. Notice that $[H^0_{\mathfrak{M}}(G)]_t$ can be written as $[H^0_{\mathfrak{M}}(G)]_t = I_t/\mathfrak{m}^{t+1}$ for some $\mathfrak{m}$-primary ideal $I_t$ such that $\mathfrak{m}^{t+1} \subseteq I_t \subseteq \mathfrak{m}^t$. Then we can write $[\overline{G}]_t = \mathfrak{m}^t/I_t$.

Taking a system of elements $f_1, \ldots, f_r$ in $\mathfrak{m}^t$ whose images in $[\overline{G}]_t$ form a $k$-basis of $[\overline{G}]_t$, we get that $a^q f_1 + \mathfrak{m}^{[q]}\mathfrak{m}^{t+1}, \ldots, a^q f_r + \mathfrak{m}^{[q]}\mathfrak{m}^{t+1}$ in $\mathfrak{m}^{[q]}\mathfrak{m}^t/\mathfrak{m}^{[q]}\mathfrak{m}^{t+1}$ are linearly independent over $k = A/\mathfrak{m}$. Actually, suppose that there exists a relation as follows:

$$\sum_{i=1}^{r} b_i(a^q f_i) \in \mathfrak{m}^{[q]}\mathfrak{m}^{t+1} \quad \text{for some } b_i \in A$$

If we put $f := \sum_{i=1}^{r} b_i f_i \in \mathfrak{m}^t$, then $a^q f \in \mathfrak{m}^{q+t+1} \subseteq I_{q+t}$. Because $a + I_1$ is $\overline{G}$-regular, this implies that $f \in I_t$. By the choice of the elements $f_1, \ldots, f_r$, we get $b_i \in \mathfrak{m}$ for all $i$ as required.

Next, using the above claim, we will complete the proof of this lemma. In order to do that, we rewrite $f(q)$ as follows:

$$f(q) = \sum_{n=0}^{q-1} \left\{ l_A(\mathfrak{m}^{[q]}/\mathfrak{m}^{[q]}\mathfrak{m}^n) - l_A(A/\mathfrak{m}^n) \right\}$$

(4.10.1)
$$= \sum_{n=1}^{q-1} \sum_{t=0}^{n-1} \left\{ \mu_A(\mathfrak{m}^{[q]}\mathfrak{m}^t) - l_A(\mathfrak{m}^t/\mathfrak{m}^{t+1}) \right\}$$

$$= \sum_{n=1}^{q-1} \sum_{t=0}^{n-1} \left\{ \mu_A(\mathfrak{m}^{[q]}\mathfrak{m}^t) - l_A([\overline{G}]_t) \right\} - \sum_{n=1}^{q-1} \sum_{t=0}^{n-1} l_A([H^0_{\mathfrak{M}}(G)]_t).$$

Since $H^0_{\mathfrak{M}}(G)$ is a module of finite length, the second term in Eq.(4.10.1) is a polynomial of $q$ with at most degree 1. Thus the required assertion easily follows from the above claim. $\square$

*Proof of Theorem* 4.8. In this proof, we also use the same notation as in the proof of Lemma 4.10. Moreover, we set $R' = R'(\mathfrak{m})$ and $\mathfrak{N} = (\mathfrak{m}, \mathfrak{m}t, t^{-1})R' = (\mathfrak{m}t, t^{-1})R'$. Then the homogeneous part of $\mathfrak{N}^{[q]}$ with degree $n$ is given as follows:

$$\left[\mathfrak{N}^{[q]}\right]_n = \begin{cases} A, & (n \leq -q), \\ \mathfrak{m}^{n+q}, & (-q+1 \leq n \leq -1), \\ \mathfrak{m}^{[q]} + \mathfrak{m}^{n+q}, & (0 \leq n \leq q-1), \\ \mathfrak{m}^{[q]}\mathfrak{m}^{n-q} + \mathfrak{m}^{n+q}, & (n \geq q). \end{cases}$$

Thus we get

$$l_{R'}(R'/\mathfrak{N}^{[q]})$$
$$= \sum_{n=-q+1}^{-1} l_A(A/\mathfrak{m}^{n+q}) + \sum_{n=0}^{q-1} l_A\left(\frac{\mathfrak{m}^n}{\mathfrak{m}^{[q]}+\mathfrak{m}^{n+q}}\right) + \sum_{n=q}^{\infty} l_A\left(\frac{\mathfrak{m}^n}{\mathfrak{m}^{[q]}\mathfrak{m}^{n-q}+\mathfrak{m}^{n+q}}\right)$$
$$= \sum_{n=0}^{q-1} l_A\left(\frac{A}{\mathfrak{m}^{[q]}+\mathfrak{m}^{n+q}}\right) + \sum_{n=q}^{\infty} l_A\left(\frac{\mathfrak{m}^n}{\mathfrak{m}^{[q]}\mathfrak{m}^{n-q}+\mathfrak{m}^{n+q}}\right)$$
$$\leq q \cdot l_A\left(\frac{A}{\mathfrak{m}^{[q]}}\right) + \sum_{n=q}^{\infty} l_A\left(\frac{\mathfrak{m}^n}{\mathfrak{m}^{[q]}\mathfrak{m}^{n-q}}\right).$$

On the other hand, since we have

$$l_R(R/\mathfrak{M}^{[q]}) = \sum_{n=0}^{q-1} l_A\left(\frac{\mathfrak{m}^n}{\mathfrak{m}^{[q]}\mathfrak{m}^n}\right) + \sum_{n=q}^{\infty} l_A\left(\frac{\mathfrak{m}^n}{\mathfrak{m}^{[q]}\mathfrak{m}^{n-q}}\right)$$

for all $q = p^e$, we obtain the required inequality by Lemma 4.10. □

The following question is natural.

**Question 4.11.** *Does an inequality* $e_{HK}(G(\mathfrak{m})) \leq e_{HK}(R(\mathfrak{m}))$ *hold? Moreover, how about* $e(A) \leq e_{HK}(R(\mathfrak{m}))$?

In case of two-dimensional Cohen–Macaulay local rings, the above question has an affirmative answer. But it is not true in higher dimension case; see Corollary 5.4.

**Proposition 4.12.** *Let* $(A, \mathfrak{m})$ *be a two-dimensional Cohen–Macaulay local ring, and let* $I \subseteq A$ *be an* $\mathfrak{m}$*-primary ideal. Then*

$$e_{HK}(R(I)) \geq \frac{e(R(I))}{2}.$$

*In particular,* $e_{HK}(R(\mathfrak{m})) \geq e(A)$.

*Proof.* Set $R = R(I)$ and $\mathfrak{M} = (\mathfrak{m}, It)$. Then $[\mathfrak{M}^{[q]}]_n = \mathfrak{m}^{[q]}I^n + I^{[q]}I^{n-q}$ for $0 \leq n < q$. Also, let $e_i(\mathfrak{m}|I)$ denote mixed multiplicities of $\mathfrak{m}$ and $I$; see e.g. [22] for details. Then we have

$$l_R(R/\mathfrak{M}^{[q]}) \geq \sum_{n=0}^{q-1} l_A(I^n/\mathfrak{m}^{[q]}I^n) \geq \sum_{n=0}^{q-1} l_A(I^n/\mathfrak{m}^q I^n)$$
$$= \sum_{n=0}^{q-1} l_A(A/\mathfrak{m}^q I^n) - \sum_{n=0}^{q-1} l_A(A/I^n)$$
$$= \sum_{n=0}^{q-1} \frac{1}{2}\left\{e_0(\mathfrak{m}|I)q^2 + 2e_1(\mathfrak{m}|I)qn + e_2(\mathfrak{m}|I)n^2\right\} - \sum_{n=0}^{q-1} \frac{e(I)}{2}n^2 + O(q^2)$$
$$= \frac{e_0(\mathfrak{m}|I) + e_1(\mathfrak{m}|I)}{2}q^3 + O(q^2)$$

for all $q = p^e$. Hence $2 \cdot e_{HK}(R(I)) \geq e_0(\mathfrak{m}|I) + e_1(\mathfrak{m}|I) = e(R(I))$. Moreover, as $e_1(\mathfrak{m}|I) \geq e(A)$ for $i = 0, 1$, we also obtain the last assertion. $\square$

**Question 4.13.** *Let $(A, \mathfrak{m})$ be a two-dimensional local ring and $I$ an $\mathfrak{m}$-primary ideal. Then does an inequality $e_{HK}(R(I)) > \frac{e(R(I))}{2}$ hold ?*

## 5. SOME EXAMPLES

In this section, we collect some results on the Hilbert-Kunz multiplicity of Rees algebras. The main aim of this section, we give a formula of $e_{HK}(R(\mathfrak{m}))$ for Veronese subring $A = k[x_1, \ldots, x_d]^{(c)}$ using Theorem 5.1.

Also, we calculate $e_{HK}(R(I))$ and $e_{HK}(R'(I))$ for a complete intersection ideal $I = (x^m, y^n)$ in $A = k[x, y]$ using Gröbner basis. Furthermore, utilizing the fact that the Hilbert-Kunz multiplicity can be described as a sort of "volume", we propose a necessary and sufficient condition for which $e_{HK}(A) = e_{HK}(R'(\mathfrak{m}))$ holds in case of two-dimensional semigroup rings.

### 5.1. Rees algberas over Veronese subrings.

Now let $A = \oplus_{n \geq 0} A_n$ be a graded ring over a field $k = A_0$. Put $\mathfrak{m} := A_+ = \oplus_{n \geq 1} A_n$. Further assume that $A = k[A_1]$. Let $\{x_1, \ldots, x_v\}$ be a $k$-basis of $A_1$ and fix it. Let $c$ be any positive integer. Let $A_1^{[cq]}$ denote the $k$-subspace generated by $x_1^{cq}, \ldots, x_v^{cq}$ of $A_{cq}$, and put $\mathfrak{m}^{[cq]} = (x_1^{cq}, \ldots, x_v^{cq})$. Note that if $c \geq 2$, then it does not hold $(x_1^{cq}, \ldots, x_v^{cq}) = (x^{cq} \,|\, x \in \mathfrak{m})$ in general. Set

(5.0.1) $\qquad \beta_n := l_A(\mathfrak{m}^n/\mathfrak{m}^{n+1}) \quad \text{and} \quad \beta_{n,cq} := l_A(\mathfrak{m}^n/\mathfrak{m}^n \cap \mathfrak{m}^{cq} + \mathfrak{m}^{n+1})$

for all $n, q \geq 0$. Then

$$\beta_n = \dim_k A_n \quad \text{and} \quad \beta_{n,cq} = \begin{cases} \beta_n & (0 \leq n \leq cq - 1) \\ \dim_k A_n/A_1^{[cq]} A_{n-cq} & (n \geq cq) \end{cases}$$

since $\mathfrak{m}^n \cap \mathfrak{m}^{[cq]} = \mathfrak{m}^{n-cq} \mathfrak{m}^{[cq]} = A_1^{[cq]} A_{n-cq} \oplus A_1^{[cq]} A_{n-cq+1} \oplus \cdots$. Note that $\beta_{n,cq} = 0$ for all $n \geq vcq$ since $\mathfrak{m}^{vcq} \subseteq (x_1^{cq}, \ldots, x_v^{cq})$ in that case.

The following theorem is a main tool for calculus of the Hilbert-Kunz multiplicity of the Rees algebra $R(\mathfrak{m})$ over Veronese subrings.

**Theorem 5.1.** *Let $A = k[A_1]$ be a homogeneous $k$-algebra with $d = \dim A \geq 2$. Let $c$ be any positive integer. Put $\mathfrak{m} = A_+$, and let $\{x_1, \ldots, x_v\}$ be a $k$-basis of $A_1$. For such a fixed system $\underline{x} = \{x_1, \ldots, x_v\}$, we define $\beta_{n,cq}$ as Eq.(5.0.1). Also, for any integers $a, k \geq 0$, we consider the following limits:*

$$I_k(a) = \lim_{q \to \infty} \frac{1}{(cq)^{d+k}} \sum_{n=0}^{aq-1} n^k \beta_{n,cq},$$

$$I_k(\infty) = \lim_{a \to \infty} I_k(a).$$

*Further, we assume that the following condition:*

(#) *The generalized Hilbert-Kunz multiplicity with respect to* $\{\underline{x}, \underline{x}t\}$

$$\lim_{n\to\infty} \frac{l_{R(\mathfrak{m})}(R(\mathfrak{m})/(x_1^n,\ldots,x_v^n,(x_1t)^n,\ldots,(x_vt)^n))}{n^{d+1}}$$

*exists; see also* [Co].

Then we have

(5.1.1) $\quad \mathrm{e}_{\mathrm{HK}}(R(\mathfrak{m})) = \dfrac{e(A) \cdot 2^{d+1}}{(d+1)!} + I_1(\infty) - 2 \cdot I_0(2c) + I_1(2c).$

*Proof.* We can regard $R(\mathfrak{m}) = A[\underline{x}t]$ as a bigraded ring $R = k[A_1t_1, A_1t_2]$ where $\deg(t_1) = \deg(t_2) = 1$. Put $\mathfrak{M} = (A_1t_1, A_1t_2)$ and

$$\mathfrak{M}^{[cq]} = (x_1^{cq}t_1^{cq}, \ldots, x_v^{cq}t_1^{cq}, x_2^{cq}t_2^{cq}, \ldots, x_v^{cq}t_2^{cq}).$$

Then the $n$th graded piece $[\mathfrak{M}^{[cq]}]_n$ of $\mathfrak{M}^{[cq]}$ (with respect to total grading) can be written as follows:

$$[\mathfrak{M}^{[cq]}]_n = \bigoplus_{\substack{i_1, i_2 \geq 0, \\ i_1+i_2=n, \\ i_1 \geq cq \text{ or } i_2 \geq cq}} A_1^{[cq]} A_{n-cq} t_1^{i_1} t_2^{i_2}.$$

Thus we have

$$l_R(R/\mathfrak{M}^{[cq]}) = \sum_{n=0}^{2(cq-1)} \dim_k A_n \times \# \left\{ (i_1,i_2) \in \mathbb{Z}^2 \;\middle|\; \begin{array}{c} 0 \leq i_1, i_2 \leq cq-1, \\ i_1 + i_2 = n \end{array} \right\}$$

(5.1.2) $\quad + \displaystyle\sum_{n=0}^{\infty} \dim_k A_n/A_1^{[cq]}A_{n-cq} \times \# \left\{ (i_1,i_2) \in \mathbb{Z}^2 \;\middle|\; \begin{array}{c} i_1, i_2 \geq 0, \\ i_1+i_2 = n, \\ i_1 \geq cq \text{ or } i_2 \geq cq \end{array} \right\}$

$$= \sum_{n=0}^{2(cq-1)} \alpha_{2,n,cq}\beta_n + \sum_{n=0}^{\infty} \alpha_{2,n}\beta_{n,cq} - \sum_{n=0}^{2(cq-1)} \alpha_{2,n,cq}\beta_{n,cq},$$

where $\alpha_{2,n} = \max\{n+1, 0\}$ and $\alpha_{2,n,cq} = \alpha_{2,n} - 2\alpha_{2,n-cq} + \alpha_{2,n-2cq}$ for all integers $n$, $q$; see also Section 3. In particular, we have

$$\alpha_{2,n,cq} = \begin{cases} n+1, & (0 \leq n \leq cq-1), \\ 2cq - n - 1, & (cq \leq n \leq 2cq-1), \\ 0, & (\text{otherwise}). \end{cases}$$

By the assumption (#), we have

$$\mathrm{e}_{\mathrm{HK}}(R(\mathfrak{m})) = \lim_{q\to\infty} \frac{1}{(cq)^{d+1}} l_R(R/\mathfrak{M}^{[cq]}).$$

From now on, we investigate each term in Eq. (5.1.2). In the proof of Lemma 3.8, if we put $c = 2$ and replace $\alpha_{d,n}$ with $\beta_n = \frac{e(A)}{(d-1)!}n^{d-1} + O(n^{d-2})$, then we get

$$(5.1.3) \quad \lim_{q \to \infty} \frac{1}{(cq)^{d+1}} \sum_{n=0}^{2(cq-1)} \alpha_{2,n,cq} \beta_n = \frac{e(A) \cdot 2!}{(2+d-1)!} S(2+d-1, 2) = \frac{e(A) \cdot 2(2^d - 1)}{(d+1)!}.$$

Since $\sum_{n=0}^{\infty} \beta_{n,cq} = l_A(A/\mathfrak{m}^{[cq]}) = O(q^d)$, we also have

$$(5.1.4) \quad \lim_{q \to \infty} \frac{1}{(cq)^{d+1}} \sum_{n=0}^{\infty} \alpha_{2,n} \beta_{n,cq} = \lim_{q \to \infty} \frac{1}{(cq)^{d+1}} \sum_{n=0}^{\infty} n \beta_{n,cq} = I_1(\infty).$$

In order to complete the proof, it is enough to show the following equality:

$$(5.1.5) \quad \lim_{q \to \infty} \frac{1}{(cq)^{d+1}} \sum_{n=0}^{2(cq-1)} \alpha_{2,n,cq} \beta_{n,cq} = -\frac{2 \cdot e(A)}{(d+1)!} + 2 \cdot I_0(2c) - I_1(2c).$$

In fact, since $\beta_{n,cq} = \beta_n$ for all $n \leq cq - 1$, we have

$$\sum_{n=0}^{2(cq-1)} \alpha_{2,n,cq} \beta_{n,cq} = \sum_{n=0}^{cq-1} (n+1) \beta_n + \sum_{n=q}^{2(cq-1)} (2cq - n - 1) \beta_{n,cq}.$$

Hence

$$\lim_{q \to \infty} \frac{1}{(cq)^{d+1}} \sum_{n=0}^{2(cq-1)} \alpha_{2,n,cq} \beta_{n,cq}$$

$$= \lim_{q \to \infty} \frac{1}{(cq)^{d+1}} \sum_{n=0}^{q-1} n \beta_n + \lim_{q \to \infty} \frac{1}{(cq)^{d+1}} \sum_{n=cq}^{2(cq-1)} (2cq - n) \beta_{n,cq}$$

$$= \lim_{q \to \infty} \frac{1}{(cq)^{d+1}} \sum_{n=0}^{cq-1} (2n - 2cq) \beta_n + \lim_{q \to \infty} \frac{1}{(cq)^{d+1}} \sum_{n=0}^{2(cq-1)} (2cq - n) \beta_{n,cq}$$

$$= 2 \lim_{q \to \infty} \frac{1}{(cq)^{d+1}} \sum_{n=0}^{cq-1} (n - cq) \beta_n + 2 \cdot I_0(2c) - I_1(2c).$$

Also, since

$$\lim_{q \to \infty} \frac{1}{(cq)^{d+1}} \sum_{n=0}^{cq-1} (n - cq) \beta_n = \lim_{q \to \infty} \frac{1}{(cq)^{d+1}} \sum_{n=0}^{cq-1} \frac{(n - cq) e(A) n^{d-1}}{(d-1)!} = -\frac{e(A)}{(d+1)!},$$

we get

$$\lim_{q \to \infty} \frac{1}{(cq)^{d+1}} \sum_{n=0}^{2(cq-1)} \alpha_{2,n,cq} \beta_{n,cq} = -\frac{2e(A)}{(d+1)!} + 2 \cdot I_0(2c) - I_1(2c).$$

Summarizing the above equalities (5.1.3), (5.1.4), and (5.1.5), we obtain the required equality. □

We are now ready to state the main theorem in this section, which is a generalization of Example 4.7. Now we will compute the Hilbert-Kunz multiplicity of the Rees algebra over Veronese subrings. Note that $A = k[x_1, \ldots, x_d]^{(c)}$ can be regarded as a homogeneous $k$-algebra with $\deg(x_1^{j_1} \cdots x_d^{j_d}) = \frac{1}{c} \sum_{i=1}^d j_i$. Putting $\underline{x} = \{x_1^{j_1} \cdots x_d^{j_d} \mid \sum j_i = c, j_1, \ldots, j_d \geq 0\}$ and $\mathfrak{m} = (\underline{x})A$, we want to apply Theorem 5.1 in this case. In Theorem 5.1, we do not know whether we need assume the condition (#) or not; see also [BC, Remark 2]. However, for affine semigroup rings, the generalized Hilbert-Kunz multiplicity always exists as is shown implicitly in the proof of [6, Theorem 2.2].

**Theorem 5.2.** *Let $c, d$ be positive integers with $c \geq d \geq 2$. Let $A = k[x_1, \ldots, x_d]^{(c)}$ be the Veronese subring where $k$ is a field of characteristic $p > 0$. Then*

$$e_{HK}(R(\mathfrak{m})) = \frac{2^{d+1} \cdot c^{d-1}}{(d+1)!} - \frac{2c - d(d-1)}{c(d+1)!} \prod_{i=1}^{d-1}(c+i).$$

*Remark 5.3.* Let $A = k[x_1, \ldots, x_d]^{(c)}$ be the Veronese subring where $k$ is a field of characteristic $p > 0$ and $c, d$ are positive integers. If we do not assume that $c \geq d$, in general, then we have

$$e_{HK}(R(\mathfrak{m})) = \frac{c^{d-1} \cdot 2^{d+1}}{(d+1)!} + \frac{1}{2c^2}\left((d+2c)\alpha_{d+1,c-1} - 2\alpha_{d+2,c-1}\right)$$
$$- \frac{2(c-1)}{cd!} \sum_{l=0}^{c-1} \alpha_{d,l} \sum_{i=0}^{N_l}(-1)^i \binom{d}{i}(2c-l-i)^d$$
$$- \frac{1}{c^2(d+1)!} \sum_{l=0}^{c-1} \alpha_{d,l} \sum_{i=0}^{N_l}(-1)^i \binom{d}{i}(2c-l-i)^{d+1},$$

where $N_l = \min\{d, 2c-l-1\}$.

*Proof of Theorem 5.2.* For simplicity, we put

$$\mathfrak{m}^{[cq]} = (x_1^{j_1 cq} \cdots x_d^{j_d cq} \mid \sum j_i = c, j_1, \ldots, j_d \geq 0)A.$$

We show that $\beta_{n,cq} = \dim_k [A/\mathfrak{m}^{[cq]}]_n$ can be written as follows:

(5.2.1) $$\beta_{n,cq} = \sum_{l=0}^{c-1} \alpha_{d,l} \sum_{i=0}^{d}(-1)^i \binom{d}{i}\alpha_{d,c(n-lq-iq)},$$

where $\alpha_{d,n} = \binom{n+d-1}{d-1}$. Note that $\beta_{n,cq} = \beta_n = \alpha_{d,cn}$ for every $n < cq$.

To see Eq. (5.2.1), for all $i_1, \ldots, i_d \geq 0$, we put
$$A_{q;i_1,\ldots,i_d} = \bigoplus_{\substack{cqi_l \leq j_l < cq(i_l+1) \\ c|j_1+\cdots+j_d}} k\ x_1^{j_1} x_2^{j_2} \cdots x_d^{j_d}.$$

Then $A = \bigoplus_{i_1,\ldots,i_d \geq 0} A_{q;i_1,\ldots,i_d}$ and $A_{q;i_1,\ldots,i_d} = x_1^{cqi_1} \cdots x_d^{cqi_d} A_{q;0,\ldots,0}$. Also, since $\mathfrak{m}^{[cq]} = \bigoplus_{i_1+\cdots+i_d \geq c} A_{q;i_1,\ldots,i_d}$, we have

$$\beta_{n,cq} = \dim_k \left[ A/\mathfrak{m}^{[cq]} \right]_n = \sum_{i_1+\cdots+i_d < c} \dim_k(A_{q;i_1,\ldots,i_d})_n$$

$$= \sum_{l=0}^{c-1} \sum_{i_1+\cdots+i_d=l} \dim_k(A_{q;i_1,\ldots,i_d})_n$$

$$= \sum_{l=0}^{c-1} \sum_{i_1+\cdots+i_d=l} \dim_k(A_{q;0,\ldots,0})_{n-lq}$$

$$= \sum_{l=0}^{c-1} \alpha_{d,l} \alpha_{d,c(n-lq),cq} \quad (\text{cf. Eq.}(3.6.1))$$

$$= \sum_{l=0}^{c-1} \alpha_{d,l} \sum_{i=0}^{d} (-1)^i \binom{d}{i} \alpha_{d,c(n-lq-iq)},$$

where the last equality follows from Eq. (3.6.2).

The following target is to calculate the following values:
$$I_0(a) = \lim_{q \to \infty} \frac{1}{(cq)^d} \sum_{n=0}^{aq} \beta_{n,cq}, \qquad I_1(a) = \lim_{q \to \infty} \frac{1}{(cq)^{d+1}} \sum_{n=0}^{aq} n\beta_{n,cq}.$$

Actually, since $\alpha_{d,cN} = \frac{c^{d-1}}{(d-1)!} N^{d-1} + O(N^{d-2})$, we have

$$I_0(a) = \sum_{l=0}^{c-1} \alpha_{d,l} \sum_{i=0}^{\min\{d,a-l\}} (-1)^i \binom{d}{i} \frac{1}{c^d} \lim_{q \to \infty} \sum_{n=lq+iq}^{aq} \frac{\alpha_{d,c(n-lq-iq)}}{q^d}$$

(5.2.2) $$= \frac{1}{cd!} \sum_{l=0}^{\min\{c-1,a\}} \alpha_{d,l} \sum_{i=0}^{\min\{d,a-l\}} (-1)^i \binom{d}{i} (a-l-i)^d.$$

Similarly, we have

$$I_1(a) = \sum_{l=0}^{c-1} \alpha_{d,l} \sum_{i=0}^{\min\{d,a-l\}} (-1)^i \binom{d}{i} \frac{1}{c^{d+1}} \lim_{q \to \infty} \sum_{n=lq+iq}^{aq} \frac{n\alpha_{d,c(n-lq-iq)}}{q^{d+1}}$$

(5.2.3) $$= \frac{1}{c^2(d-1)!} \sum_{l=0}^{\min\{c-1,a\}} \alpha_{d,l} \sum_{i=0}^{\min\{d,a-l\}} (-1)^i \binom{d}{i} \lim_{q \to \infty} \sum_{n=lq+iq}^{aq} \frac{n(n-lq-iq)^{d-1}}{q^{d+1}}$$

$$= \frac{1}{c^2(d+1)!} \sum_{l=0}^{\min\{c-1,a\}} \alpha_{d,l} \sum_{i=0}^{\min\{d,a-l\}} (-1)^i \binom{d}{i} (ad+l+i)(a-l-i)^d.$$

Now assume that $a \geq c + d$. Then $\min\{c - 1, a\} = c - 1$ and $\min\{d, a - l\} = d$ for all $l$ with $0 \leq l \leq c - 1$. Also, we recall the following property of the Stirling number of the second kind:

(5.2.4) $$\frac{1}{d!} \sum_{i=0}^{d} (-1)^{i+d} \binom{d}{i} i^k = \begin{cases} 0, & (0 \leq k < d), \\ S(d, d) = 1, & (k = d), \\ S(d+1, d) = \frac{(d+1)d}{2}, & (k = d+1). \end{cases}$$

It follows from Eq. (5.2.3) and Eq. (5.2.4) that

$$I_1(a) = \frac{1}{c^2(d+1)!} \sum_{l=0}^{c-1} \alpha_{d,l} \sum_{i=0}^{d} (-1)^{d+i} \binom{d}{i} \left[(d+1)l \cdot i^d + i^{d+1}\right]$$

$$= \frac{1}{c^2} \sum_{l=0}^{c-1} l\alpha_{d,l} \cdot S(d,d) + \frac{1}{c^2(d+1)} \sum_{l=0}^{c-1} \alpha_{d,l} \cdot S(d+1,d)$$

$$= \frac{d}{c^2} \sum_{l=1}^{c-1} \alpha_{d+1,l-1} + \frac{d}{2c^2} \sum_{l=0}^{c-1} \alpha_{d,l}$$

$$= \frac{d}{c^2} \alpha_{d+2,c-2} + \frac{d}{2c^2} \alpha_{d+1,c-1} \quad \left(\text{because } \sum_{l=0}^{N} \alpha_{d,l} = \alpha_{d+1,N}\right)$$

$$= \frac{1}{2c^2} \{(d+2c)\alpha_{d+1,c-1} - 2\alpha_{d+2,c-1}\}.$$

Similarly, if $a \geq c + d$, then Eq. (5.2.2) implies that

$$I_0(a) = \frac{1}{cd!} \sum_{l=0}^{c-1} \alpha_{d,l} \sum_{i=0}^{d} (-1)^i \binom{d}{i} (a-l-i)^d = \frac{1}{c} \sum_{l=0}^{c-1} \alpha_{d,l} S(d,d) = \frac{\alpha_{d+1,c-1}}{c}.$$

Further, assume that $c \geq d \geq 2$. Then since $2c \geq c + d$, we have

(5.2.5) $$I_1(\infty) = I_1(2c) = \frac{1}{2c^2} \{(d+2c)\alpha_{d+1,c-1} - 2\alpha_{d+2,c-1}\}, \quad I_0(2c) = \frac{\alpha_{d+1,c-1}}{c}.$$

Also, since $e(A) = c^{d-1}$, Theorem 5.1 yields that

(5.2.6) $$e_{HK}(R(\mathfrak{m})) = \frac{e(A) \cdot 2^{d+1}}{(d+1)!} + I_1(\infty) - 2 \cdot I_0(2c) + I_1(2c).$$

Substituting Eq. (5.2.5) to Eq.(5.2.6), we get

$$e_{HK}(R(\mathfrak{m})) = \frac{c^{d-1} \cdot 2^{d+1}}{(d+1)!} + \frac{1}{c^2} \{(d+2c)\alpha_{d+1,c-1} - 2\alpha_{d+2,c-1}\} - \frac{2\alpha_{d+1,c-1}}{c}$$

$$= \frac{c^{d-1} \cdot 2^{d+1}}{(d+1)!} + \frac{1}{c^2} \{d\alpha_{d+1,c-1} - 2\alpha_{d+2,c-1}\}$$

$$= \frac{c^{d-1} \cdot 2^{d+1}}{(d+1)!} + \frac{1}{c^2} \left\{d\binom{d+c-1}{d} - 2\binom{d+c}{d+1}\right\}$$

$$= \frac{2^{d+1} \cdot c^{d-1}}{(d+1)!} - \frac{2c - d(d-1)}{c(d+1)!} \prod_{i=1}^{d-1}(c+i),$$

as required. □

The following corollary gives a negative answer to Question 4.11 in case of $\dim A \geq 3$.

**Corollary 5.4.** *Let $A = k[x_1, \ldots, x_d]^{(c)}$ be the Veronese subring where $k$ is a field of characteristic $p > 0$ and $c, d$ be positive integers. Let $\mathfrak{m}$ be the unique homogeneous maximal ideal of $A$. If $d \geq 3$ and $c \geq \frac{d(d-1)}{2}$, then we have*

$$e_{HK}(R(\mathfrak{m})) \leq \frac{2^{d+1}}{(d+1)!}c^{d-1} < c^{d-1} = e(A).$$

*In fact, if we fix an integer $d \geq 3$, then $\lim_{c \to \infty} \frac{e_{HK}(R(\mathfrak{m}))}{e(A)} = \frac{2^{d+1} - 2}{(d+1)!}$.*

**Discussion 5.5.** *Under the same notation as in Theorem 5.2, we define a real continuous function $f_c : \mathbb{R} \to \mathbb{R}$ as follows:*

$$f_c(t) = \begin{cases} \dfrac{1}{c(d-1)!} \displaystyle\sum_{l=0}^{\min\{c-1,[t]\}} \alpha_{d,l} \sum_{i=0}^{\min\{d,[t]-l\}} (-1)^i \binom{d}{i}(t-l-i)^{d-1}, & \text{if } t \geq 0, \\ 0, & \text{otherwise.} \end{cases}$$

*Then $f_c(t) = 0$ for all $t \geq c + d - 1$ and one can easily check the following equality for all $a \geq 0$:*

$$I_k(a) := \lim_{q \to \infty} \frac{1}{(cq)^{d+k}} \sum_{n=0}^{aq-1} n^k \beta_{n,cq} = \frac{1}{c^k} \int_0^a t^k f_c(t)\, dt.$$

*In particular, the Hilbert-Kunz multiplicities can be represented in terms of integrals:*

$$e_{HK}(A) = \int_0^\infty f_c(t)\, dt,$$

$$e_{HK}(R(\mathfrak{m})) = \frac{c^{d-1} \cdot 2^{d+1}}{(d+1)!} + \frac{1}{c}\int_0^\infty tg(t)\, dt - 2\int_0^{2c} g(t)\, dt + \frac{1}{c}\int_0^{2c} tg(t)\, dt.$$

**5.2. Rees algberas of complete intersection ideals.**

The following example will be useful to construct counterexamples to several questions. See also Example 2.4 and Remark 4.9.

**Example 5.6.** Let $A = k[x, y]$ be a polynomial ring with two-valuables over a field $k$. Set $I = (x^m, y^n)$, where $m \geq n \geq 1$. Then we have
 (1) $e(R(I)) = n + 1$.
 (2) $e_{HK}(R(I)) = n + 1 - \dfrac{n}{m} + \dfrac{n}{3m^2}$.
 (3) $e(R'(I)) = n + 2$ $(n \geq 2)$, $= 2$ (otherwise).
 (4) $e_{HK}(R'(I)) = n + 2 - \dfrac{n}{m} - \dfrac{1}{n}$.

*Proof.* Since $R(I)$ is a binomial hypersurface, (2) follows from Lemma 1.6. Thus we must check (4) only. Set $R' := R'(I)$. Then

$$R' = k[x, y, x^m t, y^n t, t^{-1}] \cong k[x, y, z, w, t]/(x^m - zt, y^n - wt).$$

To see $e_{HK}(R') = n + 2 - \frac{n}{m} - \frac{1}{n}$, we consider the following ideal $I_q$ in $S := k[X, Y, Z, W, T]$:

$$I_q = (X^m - ZT, Y^n - WT, X^q, Y^q, Z^q, W^q, T^q)$$

for all $q = p^e$. In order to compute the length of $S/I_q$, we consider the lexicographic order ">" such that $x > y > z > w > t$. Put $c = [q/m]$ and $d = [q/n]$. Then we can easily check that

$$G = \{x^m - zt, y^n - wt, z^q, w^q, t^q, x^{q-cm}(zt)^c, y^{q-dn}(wt)^d, (zt)^{c+1}, (wt)^{c+1}\}$$

is a Gröbner basis of $I_q$ for all large $q = p^e$. Note that $in(x^m - zt)$ and $in(y^n - wt)$ are relatively prime. By the similar argument as in the proof of [4,(3.1)], we obtain the required formula for $e_{HK}(R'(I))$. $\square$

### 5.3. Other examples.

In Section 2, we have proved $e_{HK}(A) \leq e_{HK}(R'(\mathfrak{m})) \leq e_{HK}(G(\mathfrak{m}))$ for any local ring $(A, \mathfrak{m})$. So it is natural to ask when equalities hold. For example, we calculate $e_{HK}(R'(\mathfrak{m}))$ for two-dimensional affine semigroup rings of type $(A_n)$ in Example 2.4. As a result, we obtained that $e_{HK}(A) = e_{HK}(R'(\mathfrak{m}))$ if and only if $n = 1$. The following proposition gives a slight generalization of this result.

**Proposition 5.7.** *Let $s \geq 2$ be an integer. Let $S$ denote the subsemigroup of $\mathbb{Z}^2$ generated by $(a_i, b_i)$ for $i = 0, 1, \ldots, s$ where $a_i, b_j$ are integers such that $0 = a_0 < a_1 < \cdots < a_s = a$ and $0 = b_s < b_{s-1} < \cdots < b_0 = b$. Put $A = k[S] = k[x^{a_i} y^{b_i} \mid i = 0, \ldots, s]$, the affine semigroup ring of $S$ over a field $k$ of characteristic $p > 0$.*

*Then $e_{HK}(A) \leq e_{HK}(R'(\mathfrak{m}))$ and equality holds if and only if $\frac{a_i}{a} + \frac{b_i}{b} = 1$ for all $i = 0, \ldots, s$.*

*When this is the case, we also have $e_{HK}(G(\mathfrak{m})) = e_{HK}(A)$.*

*Proof.* Take the following points in $\mathbb{R}^3$: $P_{-1} = (0, 0, -1)$, $P_i = (a_i, b_i, 1)$ and $\overline{P_i} = (a_i, b_i, 0)$ for $i = 0, \ldots, s$. Let $C$ be the cone defined by the half lines $OP_i$ for $i = -1, 0, \ldots, s$ and put $V := C - \bigcup_{i=-1}^{s}(P_i + C)$. Then $vol(V) = \delta \cdot e_{HK}(R'(\mathfrak{m}))$ where $vol$ denotes the volume and $\delta = |\mathbb{Z}^2/\mathbb{Z}S|$; see [6, Theorem 2.2, Corollary 2.3].

Similarly, letting $\overline{C}$ be the cone defined by the half lines $O\overline{P_i}$ and putting $\overline{V} := \overline{C} - \bigcup_{i=-1}^{s}(\overline{P_i} + \overline{C})$ (where $\overline{P_{-1}} = P_{-1}$), we have $vol(\overline{V}) = \delta \cdot e_{HK}(A)$.

Now let $Q = (x, y, 0) \in \mathbb{R}^2$. Then $Q \in \overline{V}$ if and only if $h(Q) = 1$, where $h(Q) = \max\{z | (x, y, z) \in V\} - \min\{z | (x, y, z) \in V\}$ since each $O\overline{P_i}$ is the projection of $OP_i$ to $xy$-plane. This implies that $vol(V) \geq vol(\overline{V})$ and $e_{HK}(R'(\mathfrak{m})) \geq e_{HK}(A)$.

Now suppose that $e_{HK}(R'(\mathfrak{m})) = e_{HK}(A)$ holds. Then $Q \notin \overline{V}$ if and only if $h(Q) = 0$ or $h(Q)$ cannot be defined. Thus two lines $P_{i-1} + tOP_s$ and $P_i + t'OP_0$ must intersect for $i = 1, \ldots, s$. That is, $(a_{i-1}, b_{i-1}, 1) + t(a, 0, 1) = (a_i, b_i, 1) + t'(0, b, 1)$. This yields $a_{i-1}/a + b_{i-1}/b = a_i/a + b_i/b$. Therefore $a_i/a + b/b_i = b/b = 1$ for each $i$.

Conversely, suppose that $a_i/a + b_i/b = 1$ holds for all $i$. Then we can regard $A$ as a homogeneous $k$-algebra with $\deg(x) = 1/a$ and $\deg(y) = 1/b$. Hence we have $e_{HK}(G(\mathfrak{m})) = e_{HK}(R'(\mathfrak{m})) = e_{HK}(A)$ by Example 2.3. □

**Corollary 5.8.** *Let $A$ be a two-dimensional normal semigroup ring. If $e_{HK}(R'(\mathfrak{m})) = e_{HK}(A)$, then it is isomorphic to a Veronese subring.*